\documentclass[11pt]{amsart}

\usepackage{url}
\usepackage{psfrag}
\usepackage{graphicx}
\usepackage{comment}
\usepackage{array}
\usepackage{amsmath,amsfonts,amssymb,amsxtra,amsthm}
\usepackage{mathrsfs}
\usepackage{lineno}
\usepackage{verbatim}
\usepackage{xy}
\usepackage{times}
\xyoption{all}
\usepackage{stmaryrd}
\usepackage{appendix}
\usepackage{xr}

\DeclareMathOperator{\Aut}{Aut}

\DeclareMathOperator{\Hom}{Hom}

\DeclareMathOperator{\Ker}{Ker} 

\DeclareMathOperator{\Mod}{Mod}

\DeclareMathOperator{\stabker}{\underrightarrow{\Ker}}

\DeclareMathOperator{\trace}{tr}

\def\jsj{\textsc{JSJ}}

\def\hnn{\textsc{HNN}}

\def\wt#1{\widetilde{#1}}

\def\sl2c{\ensuremath{{SL}(2,\mathbb{C})}}
\def\t1sl2c{{\mathfrak sl}_2\mathbb{C}}
\def\free{\ensuremath{\mathbb{F}}}

\def\onto{\twoheadrightarrow}
\def\into{\hookrightarrow}

\newcommand{\adjoin}[2]{\ensuremath{#1\!\left[#2\right]}}

\def\cent{Z}

\makeatletter
\newcommand{\Doubletwo}[3]{ \left\{ #1 \mid #3\right\} }
\newcommand{\Doubleone}[1]{ \left\{ #1 \right\} }
\newcommand{\Doublethr}[3]{ \left\{ #1 \right\}_{#3} }
\newcommand{\set}[1]{%
\@ifnextchar:{\Doubletwo{#1}}{\@ifnextchar_{\Doublethr{#1}}{\Doubleone{#1}}}%
}
\makeatother

\makeatletter
\newcommand{\grouptwo}[3]{ \langle #1 \mid #3\rangle }
\newcommand{\groupone}[1]{ \langle #1 \rangle }
\newcommand{\group}[1]{%
\@ifnextchar:{\grouptwo{#1}}{\groupone{#1}}%
}
\makeatother

\def\eval{\ensuremath{\mathrm{ev}}}
\def\repvar{\ensuremath{\mathcal{R}}}

\newcommand{\term}[1]{{\emph{#1}}}
\newcommand{\scare}[1]{`#1'}

\newtheorem{theorem}{Theorem}[section]
\newtheorem{lemma}[theorem]{Lemma}

\theoremstyle{definition}
\newtheorem{definition}[theorem]{Definition}

\theoremstyle{remark}
\newtheorem{remark}[theorem]{Remark}

\newcommand{\mnote}[1]{}

\title[Bounding strict resolutions]{Krull dimension for limit groups
  {I}:\\
{B}ounding strict resolutions}

\author{Larsen Louder}
\address{Department of Mathematics\\
University of Michigan \\
Ann Arbor, MI 48109-1043\\
USA}
\email[Larsen Louder]{llouder@umich.edu, lars@d503.net}



\date{\today}

\begin{document}

\begin{abstract}
  This is the first paper in a sequence on Krull dimension for limit
  groups, answering a question of Z.~Sela. In this paper we show that
  strict resolutions of a fixed limit group have uniformly bounded
  length. The upper bound plays two roles in our approach. First, it
  provides upper bounds for heights of analysis lattices of limit
  groups, and second, it enables the construction of
  \jsj--respecting sequences in the sequel.
\end{abstract}

\maketitle

\section{Introduction}

\thispagestyle{empty}

\par We start by drawing the analogy between limit groups and
coordinate rings of algebraic varieties. To an algebraic subset $V$ of
$\mathbb{C}^n,$ one attaches the coordinate ring. The points of $V$
are in one-to-one correspondence with ring homomorphisms
$\adjoin{\mathbb{C}}{V}\to\mathbb{C}$: If $\adjoin{\mathbb{C}}{V}$ is
the quotient of $\adjoin{\mathbb{C}}{x_i}$ by a radical ideal $I,$
then an assignment $x_i\mapsto z_i\in\mathbb{C},$ $(z_i)$ a point of
$V,$ vanishes on $I$ and gives a ring homomorphism. Conversely, one
associates to such a ring homomorphism $f$ the tuple $(f(x_i))\in V$.

\par The replacement for the notion of a radical ideal in algebraic
geometry over the free group is that of a residually free group. Going
further, we need to extend the notion of prime ideal as well. A
\term{zero-divisor} in a finitely generated group $G$ is a tuple of
elements of $G$ such that every homomorphism $G\to\free$ kills at
least one element of the tuple. A group is \term{residually free} if
it has no singleton zero divisors. A group is
\term{$\omega$--residually free}, or is a \term{limit group}, if it
has no zero divisors. It will be convenient to use a slightly
different definition of a limit group.

\begin{definition}[Limit Group \cite{sela::dgog1}]
  \label{limitgroup}
  A sequence of homomorphisms $f_n\colon G\to\Gamma$ is stable if, for
  all $g\in G,$ there exists $n_g$ such that $f_n(g)=1$ for $n>n_g$ or
  $f_n(g)\neq 1$ for $n>n_g$. The \term{stable kernel} of a stable
  sequence of homomorphisms is the set of elements which have trivial
  image for large $n,$ and is denoted $\stabker(f_n)$. The quotient of
  $G$ by the stable kernel of a stable sequence $f_n$ is a
  \term{$\Gamma$--limit group}.
\end{definition}

\par In this paper we study $\Gamma$--limit groups for
$\Gamma\cong\free$. A sequence $f_n\colon G\to\free$ \term{converges
  to} $G$ if $\stabker(f_n)=\set{1}$. If $G$ is $\omega$--residually
free then clearly there exists a sequence of homomorphisms $f_n\colon
G\to\free$ converging to $G$. That a $\free$--limit group is
$\omega$--residually free is~\cite[Theorem~4.6]{sela::dgog1}, and
follows from finite presentability.

\par The object analogous to the coordinate ring for an algebraic
subset $V$ of $\free^n$ is simply a finitely generated residually free
group $G$ which requires only $n$ generators. With little work one
establishes a one-to-one correspondence between points of $V$ and
elements of $\Hom(G,\free)$. If an algebraic set $V$ is irreducible
then $\adjoin{\mathbb{C}}{V}$ has no zero divisors. Likewise, for the
set $\Hom(G,\free)$ to be irreducible, $G$ must have no zero divisors,
hence is a limit group.

\par The \term{Krull dimension} of an algebraic set is the supremum of
lengths of chains of irreducible subvarieties, thus it is natural to
ask whether any sequence of proper epimorphisms of limit groups,
beginning with a fixed free group of rank $n,$ terminates in a uniform
number of steps. Exhibition of such a bound would establish the
finiteness of the Krull dimension of $\free^n$. In this paper we make
the first move toward showing this.

\par The key to constructing limit groups is the notion of a strict
resolution. The group $\Mod(L)<\Aut(L)$ is defined in the next
section.

\begin{definition}[Strict; Strict resolution]
  \label{strict}
  A homomorphism $\pi\colon L\to L'$ is \term{strict} if for every
  sequence of homomorphisms $f_n\colon L'\to\free$ converging to $L',$
  there exists a sequence of automorphisms $\phi_n\in\Mod(L)$ such
  that the sequence $f_n\circ\pi\circ\phi_n$ converges to $L$.

  A sequence of proper epimorphisms $L\onto L_1\dotsb\onto\free_k$ is
  a \term{resolution} of $L$.  If every homomorphism appearing in a
  resolution is strict then the resolution is a \term{strict
  resolution}.
\end{definition}

\par By~\cite[Proposition~5.10]{sela::dgog1} every limit group
admits a strict resolution. The first step in our approach to showing
that limit groups have finite Krull dimension is to show that any
strict resolution of limit groups terminates in a uniformly bounded
number of steps.

\begin{theorem}[Strict resolutions have bounded length]
  \label{strictlengthbound}
  Let $\free_n\onto L_0\onto\dotsb\onto L_k$ be a sequence of
  proper strict epimorphisms of limit groups. Then $k\leq 3n$.

  Let $\free_n\onto L_1\onto\dotsb\onto L_k\cong\free_m$ be a
  sequence of proper strict epimorphisms of limit groups. Then $k\leq
  3(n-m)$.
\end{theorem}

\par It is well known that any resolution, strict or not, has finite
length~\cite[Proposition~5.1]{sela::dgog1},~\cite[Corollary~1.9]{bf::lg}. The
importance of strict resolutions lies in the fact that they witness a
group as a limit group. See Definition~\ref{strict}. The main use of
Theorem~\ref{strictlengthbound} is along the way
to~\cite[\ref{STABLE-thr:alignmenttheorem}]{louder::stable}, which
says, roughly, that if there exist arbitrarily long sequences of
proper epimorphisms of rank $n$ limit groups, then there exist
arbitrarily long sequences of rank $n$ limit groups such that all
groups in the sequence share the same \jsj\ decomposition.

\subsection*{Acknowledgments}

The author would like to thank Mladen Bestvina, for wondering aloud if
limit groups have irreducible representation varieties. Thanks must
also go to Ben McReynolds and Chlo\'e Perin for a number of helpful
suggestions.

\section{Splittings and Lifting Automorphisms}

\par The following definition is somewhat non-standard. Its utility
lies in the fact that it streamlines the statements and proofs of the
lemmas leading up to Theorem~\ref{autspreservecomponents}.

\begin{definition}
\label{typesofsplittings}
  Fix a finitely generated group $G$. A splitting
  $G= G_1*_{E}G_2$ over a finitely generated abelian
  subgroup $E$ is an \term{amalgamation}. A splitting
  $G= G'*_{E}$ over a finitely generated abelian subgroup
  $E$ is an \term{\hnn\ extension}. In either case, $E$ is
  allowed to be trivial, as is $G_1$ in the former.
\end{definition}

An important subgroup of the automorphism group of a freely
indecomposable finitely generated group is the \term{modular group}
$\Mod$. To make our exposition as efficient as possible, we define the
modular group to be the subgroup of the automorphism group generated
by the following \term{elementary automorphisms}:

\begin{itemize}
\item Automorphisms from amalgamations: $G= G_1*_{E}G_2$: The Dehn
  twist in $e\in \cent_{G_2}(E)$ is the automorphism of $G$ which is
  the identity on $G_1$ and which is conjugation by $e$ on $G_2$. This
  automorphism is denoted $\tau_{\Delta,e}$.

\item Automorphisms from \hnn\ extensions: $G'*_{E}$: Let $t$ be the
  stable letter such that $t E t^{-1}= E_1$ and $e\in
  \cent_{G'}(E)$. Then $\tau_{\Delta,e}$ is the automorphism which is
  the identity on $G'$ and maps $t$ to $te$.
\end{itemize}

\par The modular group as defined above is the group of
\scare{geometric} automorphisms arising from one edged splittings. If
$G$ is freely indecomposable, then $\Mod(G)$ agrees with the modular
group as defined in~\cite{bf::lg} and~\cite{sela::dgog1}.

\par Given an \hnn\ extension or amalgamation of $G,$ we build a
free group $\wt{G}$ and a homomorphism
$\wt{G}\onto G$ so that the automorphism of $G$
engendered by the splitting lifts to an automorphism of
$\wt{G}$.

\begin{definition}[Lifts of Automorphisms, $\wt{G}$]
\label{liftsofautomorphisms}
  Let $\Delta$ be a splitting of $G$ as $G=G_1*_{E}G_2$ over a
  finitely generated abelian group $E,$ $e\in\cent_{G_2}(E)$. ($G_1$
  and $G_2$ are finitely generated.) Let $F_1=\langle
  x_1,\dotsc,x_{n_1}\rangle$ and $F_2=\langle
  y_1,\dotsc,y_{n_2}\rangle$ be free groups with homomorphisms
  $\pi_i\colon F_i\onto G_i$. Define $\wt{G}$ to be the group
  $F_1*F_2,$ calling this free splitting $\wt{\Delta},$ and let
  $\pi=\pi_1*\pi_2\colon\wt{G}\onto G$ be the obvious
  surjection. This particular free factorization of $\wt{G}$ is
  the lift of $\Delta$. Choose $\wt{e}\in F_2$ such that
  $\pi_2(\wt{e})=e$. Then the automorphism
  $\tau_{\wt{\Delta},\wt{e}}$ of $\wt{G}$ makes
  the following diagram commute.

\label{liftdiagram}
\[
    \xymatrix{%
      \wt{G}\ar[r]^{\tau_{\wt{\Delta},\wt{e}}}\ar[d]^{\pi} & \wt{G}\ar[d]^{\pi}\\
     G\ar[r]^{\tau_{\Delta,e}} &  G }
    \eqno{(\diamondsuit)}
\]

\par If $\Delta$ is an \hnn\ extension $G=G'*_{E}$ we define $\wt{G}$
similarly. Choose $F'=\langle x_1,\dotsc,x_n\rangle,$ a homomorphism
$\pi'\colon F'\onto G',$ and $\wt{G}$ the \hnn\ extension $F'*\langle
\wt{t}\rangle$. Set $\pi(\wt{t})=t$. Choose a lift $\wt{e}\in F'$ of
$e$ and define an automorphism $\tau_{\wt{\Delta},\wt{e}}$ which is
the identity on $F'$ and maps $\wt{t}$ to $\wt{t}\cdot\wt{e}$.

\end{definition}


\section{Dehn twisting $\repvar(\free)$}

\par Our approach to Theorem~\ref{strictlengthbound} is to analyze the
action of the modular group of a limit group on its $\sl2c$
representation variety.

\par The set of homomorphisms of a finitely generated group to a
complex algebraic Lie group is an affine subvariety of $\mathbb{C}^n$
for some $n$ \cite{cullershalen0}. The $\sl2c$ representation
variety $\Hom(G,\sl2c)$ will be denoted by $\repvar(G)$. For $g\in G,$
$\eval_g$ will denote the evaluation map $\repvar(G)\to\sl2c$.

\par We saw above that if $\tau_{\Delta,e}$ is elementary then it
lifts to an elementary automorphism of $\wt{G}$. The
representation variety $\repvar(G)$ has a natural embedding
\[\repvar(G)\into\repvar(\wt{G}) =\repvar(\free_n) 
=\sl2c^n\subset\mathbb{C}^{4n}\] where $n$ is the rank of
$\wt{G}$. The lift $\tau_{\wt{\Delta},\wt{e}}$
acts on $\repvar(\wt{G})$ and by commutativity of
($\diamondsuit$) the restriction
$\tau_{\wt{\Delta}}\vert_{\repvar(G)}$ agrees with
$\tau_{\Delta,e}$.

\par The exponential map
$\exp\colon M(2,\mathbb{C})\to GL(2,\mathbb{C})$ is
given by the formula
\[\exp(M)=\sum_{i=0}^{\infty}\frac{M^i}{i!}\] 
This power series converges everywhere. The following lemma is well
known. See~\cite[Chapter~1, Example~9]{rossmann} for the computation
in the real case.

\begin{lemma}
  The exponential map is biholomorphic in a neighborhood
  of all points $v\neq 0\in\t1sl2c\subset M(2,\mathbb{C})$ such that
  $\exp(v)\neq\pm I_2$.
\end{lemma}

Since all maps we deal with are either polynomials or are power series
in matrices of polynomials we suppress mention of the ambient space
$\mathbb{C}^k$.


\par The image of an edge group under a representation $\rho$
(usually) lives in a 1-parameter subgroup of $\sl2c$. We can
therefore twist $\rho$ by elements in the 1-parameter subgroup to
produce new representations. It turns out that the twisted
representations can be chosen to vary analytically in $\rho$ and a
parameter $z$ as long as $z$ is chosen carefully and the
representation $\rho$ doesn't map the edge group to an element whose
trace is $-2$.

\par To get the ball rolling we need to know where we can take
logarithms and how to define small pieces of 1-parameter subgroups.

\begin{lemma}
  \label{complexlemma}
  Let $P_{\epsilon}={N}_{\epsilon}(\set{x+0i}:{0\leq x\leq
  1})\subset\mathbb{C}$ be the epsilon neighborhood of
  $\left[0,1\right]\subset\mathbb{R}\subset\mathbb{C}$.

  For all $g\in\sl2c,$ $\trace(g)\neq-2,$ there is an element
  $v_g\in\t1sl2c,$ a neighborhood $U_g\subset\sl2c,$
  $v_g\in\wt{U}_g(=\log(U_g)$), and $\epsilon>0$ such
  that
  \begin{itemize}
    \item $P=P_{\epsilon}$
    \item $\exp(v_g)=g$
    \item
      $\exp(P\cdot\wt{ U}_g)\subset\sl2c\setminus\set{-I}$
    \item $\exp\vert_{N_{\epsilon}(\left[0,1\right]\cdot v_g)}$ is
      biholomorphic onto its image.
    \item $P\cdot\wt{ U}_g\subset
      N_{\epsilon}(\left[0,1\right]\cdot v_g)$
  \end{itemize}
\end{lemma}

\begin{proof}
  The argument is an easy adaptation of the fact that if $N\subset
  M$ is an embedded smooth, compact, submanifold and $f\colon M\to
  M'$ is smooth, injective on $N,$ and a local diffeomorphism at
  every point of $N,$ then $f$ is a diffeomorphism on a neighborhood
  of $N$.
\end{proof}

\begin{definition}
  \label{standardneighborhood}
  A tuple $S=(U_g,\wt{U}_g,P)$ satisfying
  Lemma~\ref{complexlemma} is a \term{standard neighborhood} of $g$ in
  $\sl2c$. Note that $S$ involves a particular choice of the logarithm
  $\log\colon U_g\to\t1sl2c$. If $h\in U_g$ then $\log(h)$ shall
  be taken to be the element $v_h\in\wt{U}_g$ such that
  $\exp(v_h)=h$.
\end{definition}

\begin{lemma}
  \label{autfnanalysis}
  Let $\langle x_1,\dotsc,x_n,y_1,\dotsc,y_m\rangle=\free,$
  $e\in\langle y_i\rangle,$ $\rho\in\Hom(\free,\sl2c)=\repvar(\free)$
  such that $\trace(\rho(e))\neq-2$. Choose a triple $(e,S,V)$ such
  that $\rho\in V$ and $\eval_e(V)\subset U_{\rho(e)}$. Then
  the map $\tau_H\colon V\times P\to\repvar(\free)$ defined by

  \[(\eta,z)\xrightarrow{\tau_H}
    \Big\lbrace \begin{array}{cll} x_i&\mapsto&\eta(x_i)\\
      y_i&\mapsto&\eta(y_i)^{\exp(z\cdot\log(\eta(e)))}
    \end{array} \]
  is holomorphic. Similarly, if $\free=\langle y_1\ldots
  y_{n-1}\rangle*\langle t\rangle=F*\langle t\rangle,$ $e\in F,$ then,
  after choosing an appropriate triple $(e,S,V),$ the map
  $\tau_H\colon V\times P\to\repvar(\free)$ defined by
  \[(\eta,z)\xrightarrow{\tau_H}\Big\lbrace\begin{array}{cll}
  x_i&\mapsto&\eta(x_i)\\
  t&\mapsto&\eta(t)\exp(z\cdot\log(\eta(e)))
  \end{array}\] is holomorphic.
\end{lemma}

\begin{proof}
  $\tau_H$ is the composition of holomorphic functions.
\end{proof}

\section{After lifting, restricting}

\par Let $\Delta$ be a one-edged splitting of $G,$
$\tau_{\Delta,e}$ an elementary automorphism, $\wt{G}$ the
lifted group and $\tau_{\wt{\Delta},\wt{e}}$ the lift of
$\tau_{\Delta,e}$ to $\wt{G}$. There is a natural inclusion
$\repvar(G)\subset\repvar(\wt{G})$. Since
$\repvar(\wt{G})\subset\sl2c^M\subset\mathbb{C}^{4M},$ if
$V\subset\repvar(G)$ is an open set, then
$V=\repvar(G)\cap W$ for some open subset
$ W\subset\mathbb{C}^{4M}$. If $\varphi\colon W\to\mathbb{C}$ is
analytic, its restriction to $V$ is analytic by definition.

\begin{lemma}
\label{maintwistlemma}
  Let $(\wt{e},S,V)$ be as in Lemma~\ref{autfnanalysis} and
  choose $\rho\in V\subset\repvar(G)$. Regard $\rho$ as a
  representation of $\wt{G}$ by inclusion,
  $S=(U_g,\wt{U}_g,P)$. Then
  $\tau_H((V\cap\repvar(G))\times P)\subset\repvar(G)$.

  Suppose $\repvar(G)_1,\dotsc,\repvar(G)_k$ are the
  irreducible components of $\repvar(G)$ containing $\rho$. If
  $V_i=V\cap\repvar(G)_i$ is irreducible as an analytic
  variety, then $\tau_H(V_i\times P)\subset\repvar(G)_i$.
\end{lemma}

\begin{proof}
  $\tau_H$ is cooked up in such a way that
  $\tau_H\vert_{V\cap\repvar(G)\times P}$ has image in
  $\repvar(G)$. We prove the lemma for elementary automorphisms
  arising from amalgamations. The argument in the case of an
  \hnn\ extension is identical.

  Let $K_i=\Ker(\wt{G}_i\onto G),$ and
  $\set{r_j^i}_{j=1..\infty}$ an enumeration of $K_i$.  Since finitely
  generated rings of polynomials over $\mathbb{C}$ are Noetherian,
  there exists $k<\infty$ such that
  \[\rho_i\in\repvar(G_i)\iff\forall j\leq k\left(\eval_{r^i_j}(\rho_i)=I_2\right)\]
  The inclusions $E\into G_i$ induce restriction maps
  $\repvar(G_i)\to\repvar(E)$.  Since $G$ is the pushout of the
  diagram $\set{E\into G_i},$ the representation variety $\repvar(G)$
  is the pullback of the diagram $\set{\repvar(G_i)\to\repvar(E)},$
  and we identify $\repvar(G)$ with the set of pairs
  $(g,h)\in\repvar(G_1)\times\repvar(G_2)$ such that the restrictions
  $g\vert_E$ and $h\vert_E$ agree. Since $E$ is finitely generated
  there are relations $g^l_i\in F_i,$ $l=1..m,$ corresponding to
  generators of $E,$ such that

  \[\rho=(\rho_1,\rho_2)\in\repvar(G)\subset\repvar(G_1)\times\repvar(G_2)\iff\forall l\left(\eval_{g_l^1}(\rho_1)=\eval_{g_l^2}(\rho_2)\right)\]

  The Dehn twists $\tau_H$ clearly preserve the relations
  $\eval_{r_i^j}$ and since $e\in\cent_{G_2}(E),$
  \[
    \rho_2(g_l^2)^{\exp(z\cdot\log(\eta(e)))}=\rho_2(g_l^2)
    \] 
  for all $g_l^2$. Since $\tau_H$ doesn't change the values of this
  finite set of equations, $\tau_H(V\times P)\subset\repvar(G)$.

  The intersection of an irreducible algebraic variety with an open
  subset of $\mathbb{C}^{4n}$ is an analytic variety. If the
  intersection $V_i=V\cap\repvar(G)_i$ is irreducible then
  $V_i\times P$ is also irreducible. The image of an irreducible
  complex analytic variety under a holomorphic map has irreducible
  closure (preimages of closed sets are closed), hence must have image
  in an irreducible component of the range, in this case
  $\repvar(G)$. Thus, since $V_i\nsubseteq\cup_{j\neq
  i}\repvar(G)_j,$ $\tau_H(V_i\times
  P)\subset\repvar(G)_i$.
\end{proof}

\begin{definition}
  Let $\repvar_2(G)$ be the union of the irreducible components of
  $\repvar(G)$ such that for all $g\in G,$ $\trace(\eval_g(\text{ }))+2$
  doesn't vanish. A component $V$ of $\repvar(G)$ is a component of
  $\repvar_2(G)$ if, for every $g\in G,$ there is $\rho\in V$ such
  that $\trace(\rho(g))\neq-2$.
\end{definition}

\begin{theorem}
\label{autspreservecomponents}
  Suppose $\tau_{\Delta,e}$ is an elementary automorphism of $G,$
  $\rho\in\repvar(G)_i,$ and $\trace(\rho(e))\neq -2$. Then
  $\tau_{\Delta,e}(\repvar(G)_i)=\repvar(G)_i$.

  The modular group acts trivially on the set of irreducible
  components of $\repvar_2(G)$.
\end{theorem}

\begin{proof}
  By Lemma~\ref{maintwistlemma} for an appropriate triple
  $(\wt{e},S,V),$ the 1-parameter family of Dehn twists $\tau_H$ maps
  $V_i\times P$ to $\repvar(G)_i$ as long as $V_i$ is irreducible. The
  restrictions $\tau_H\vert_{V_i\times\set{0}}$ and
  $\tau_H\vert_{V_i\times\set{1}}$ agree with the restrictions of
  $id_{G}$ and $\tau_{\Delta,e},$ respectively. Thus $\tau_{\Delta,e}$
  maps $V_i$ to $\repvar(G)_i,$ and since $\tau_{\Delta,e}$ is an
  automorphism, $\tau(V_i)$ cannot lie in the intersection
  $\repvar(G)_i\cap(\bigcup_{j\neq i}\repvar(G)_j)$. Thus
  $\tau(\repvar(G)_i)$ shares a point with
  $\repvar(G)_i\setminus(\bigcup_{j\neq i}\repvar(G)_j),$ hence
  $\tau(\repvar(G)_i)=\repvar(G)_i$.

  If $\tau_{\Delta,e}$ is an elementary automorphism of $G$ and
  $V$ is a component of $\repvar_2(G),$ then there is a
  representation $\rho\in V$ such that $\trace(\rho(e))\neq-2$. By
  the above, $\tau(V)=V$. Since the modular group is generated
  by elementary automorphisms, the claim holds.
\end{proof}

\par If $\phi\colon G\to H$ is a homomorphism then
$\phi^{-1}(\repvar_2(H))\subset\repvar_2(G)$. If $V$ is a
component of $\repvar_2(H)$ then $\phi^{-1}(V)$ is contained
in some irreducible component $W$ of $\repvar(G)$. If
$g\in G$ and $\phi(g)\neq 1$ then $\trace(\eval_g(\text{ }))+2$ doesn't
vanish on $W$ since it doesn't vanish on $V$. If $\phi(g)=1$
the same holds since $\trace(I_2)=2$.


\section{Application to strict resolutions of limit groups}

\par In this section we give our main application of
Theorem~\ref{autspreservecomponents}: a bound on the length of a
strict resolution of a limit group which depends only on its
rank. Limit groups possess two qualities which make the theory
developed so far useful: many maps to free groups, which have large
representation varieties, and large automorphism groups generated by
Dehn twists in one-edged splittings.

\par We give only enough definitions to make sense of the statement of
the theorem and its proof. They will be very economical. For more
information on limit groups see the exposition by Bestvina and Feighn
\cite{bf::lg} or Sela's original work \cite{sela::dgog1}.

\par A homomorphism $\free\to\sl2c$ is \term{nondegenerate} if it is
injective and every element has image with trace not equal to $-2$.

\begin{definition}
  The \term{essential subvariety}, $\repvar_{e}(L),$ comprises the
  irreducible components $V$ of $\repvar(L)$ such that there is a
  nondegenerate $i\in\repvar(\free)$ and a sequence $f_n\colon
  L\to\free$ converging to $L,$ such that $i\circ f_n\in V$.
\end{definition}

\par The important feature the essential subvariety has is that for
all $g\in L$ the evaluation map $\eval_g$ takes non-identity values
on every component. The essential subvariety for a limit group is
non-empty since $\free$ has a nondegenerate embedding in $\sl2c,$
limit groups are $\omega$--residually free, and $\repvar(L)$ has
finitely many components. Since there exist nondegenerate elements of
$\repvar(\free),$ $\repvar_{e}(L)\neq\emptyset$. By definition
$\repvar_{e}(L)\subset\repvar_2(L)$.

\begin{lemma}
  \label{basicdimensioninequality}
  If $\pi\colon L\onto L'$ is strict then
  $\repvar_{e}(L')\subset\repvar_{e}(L)$ and if
  $V'\subset V$ are irreducible components of
  $\repvar_{e}(L')$ and $\repvar_{e}(L),$ respectively,
  then $\dim V'<\dim V$. Thus
  $\dim\repvar_{e}(L')<\dim\repvar_{e}(L)$.
\end{lemma}

\begin{proof}
  To prove the first part of the claim, choose an irreducible
  component $V'$ of $\repvar_{e}(L')$ and an irreducible
  component $V$ of $\repvar(L)$ containing $V'$. We show
  that $V\subset\repvar_{e}(L)$.

  Choose a sequence $f_n\in\Hom(L',\free)$ converging to $L',$ and a
  nondegenerate $i\in\repvar(\free)$ such that $i\circ f_n$ is
  contained $V'$ for all $n$. Choose $\phi_n\in\Mod(L)$ such that
  $g_n=f_n\circ\pi\circ\phi_n$ converges to $L$. Since
  $\repvar_2(L')\subset\repvar_2(L)$ we have
  $\repvar_{e}(L')\subset\repvar_2(L)$. By
  Theorem~\ref{autspreservecomponents}, the $\phi_n$ fix $V$ setwise
  and we have $i g_n\in V$ for all $n,$ hence $V$ is a component of
  $\repvar_{e}(L)$ and we have the claim.

  To prove the second claim choose some $g\neq
  1\in\Ker(L\to L')$. By the definition of the essential
  subvariety, the evaluation map $\eval_g$ doesn't vanish on any
  irreducible component of $\repvar_{e}(L)$ and, since
  $g\in\Ker(\pi),$ $\eval_g$ vanishes on $\repvar_{e}(L'),$ hence the
  inequality on dimensions of irreducible components is strict.

  The last statement follows immediately from the first two.
\end{proof}

\begin{proof}[Proof of Theorem~\ref{strictlengthbound}]
  We prove only the second assertion. Let
  $\free_n\onto L_1\onto\dotsc\onto L_k\cong\free_m$ be a strict
  resolution and consider the sequence of essential subvarieties
  \[\repvar_{e}(L_k)\cong\sl2c^m\subset\dotsb\subset\repvar_{e}(L_1)\subset\sl2c^n\]
  Let $d_i=\dim(\repvar_{e}(L_i))$. Since
  $L_i\onto L_{i+1}$ is strict and proper, by
  Lemma~\ref{basicdimensioninequality}, $d_i>d_{i-1}$. Since $d_k=3m$
  and $d_1\leq 3n,$ $k$ must be at most $3(n-m)$.
\end{proof}

\section{Accessibility for limit groups}

\par Sela proves~\cite[Theorem~4.1 and Proposition~4.3]{sela::dgog1}
that the height of the cyclic analysis lattice of a limit group is
quadratic in the first Betti number. One application of
Theorem~\ref{strictlengthbound} is the following companion to Sela's
result. This kind of argument is used
in~\cite[Theorem~\ref{STABLE-lem:depthbound}]{louder::stable}
and~\cite[Theorem~\ref{ROOTS-analysislatticebound}]{louder::roots}.

\begin{definition}[Analysis lattice]
  The abelian (cyclic) analysis lattice of a finitely generated group
  $G$ is the following tree of groups. Free and abelian groups
  have no children.
  \begin{itemize}
    \item Level 0 consists of $G$.
    \item Level 0.1 consists of the freely indecomposable free factors
      of $G$ and the free group of some Grushko free factorization
      of $G$. The groups in level 0.1 are the children of the node
      labeled $G$.
    \item Level 1 consists of the rigid, abelian, and quadratically
      hanging subgroups in the abelian (cyclic) \jsj's of the freely
      indecomposable free factors of $G$ at level 0.1. The parent of a
      group at level 1 is the freely indecomposable free factor it is
      a subgroup of.
    \item Level 1.1 is constructed exactly as level 0.1 was.
\end{itemize}
\end{definition}

\begin{theorem}
  The height of the abelian analysis lattice of a limit group is
  bounded by three times its rank.
\end{theorem}

\begin{remark}
  Before we begin, note if that $L_v$ is a vertex group of the abelian
  \jsj\ decomposition of $L$ then the restriction of every modular
  automorphism to $L_v$ agrees with the restriction of an inner
  automorphism since $ L_v$ is elliptic in every splitting of $L$ over
  an abelian subgroup. If $\pi\colon L\to L'$ is strict and
  $g\in\ker(L_v\to L'),$ then $g$ is in the kernel of every element of
  $\pi\circ\Mod(L)$. Since $\pi$ is strict $g$ must be the identity
  element, therefore $L_v$ embeds in $L'$.
\end{remark}

\begin{proof}
  Let $L$ be a limit group generated by $n$ elements. We prove
  something slightly more general: that the height of the abelian
  analysis lattice is bounded linearly by the length of the shortest
  strict resolution $L$ admits. Observe that if $L'< L$ then any
  strict resolution of $L$ restricts to a strict resolution of $L'$.

  Let $L=L_1\onto L_2\onto\dotsb\onto L_k$ be a shortest strict
  resolution of $ L$. If $k=1$ then $L$ has height $0$ and we may
  stop. Otherwise, let $L''$ be a group at level 1 of the analysis
  lattice. If $L''$ is free or abelian it has height 0 so we may stop
  the procedure. If not, then by the previous paragraph $L''$ embeds
  in $L_2,$ and by induction the height of the abelian analysis
  lattice of $L_2$ is at most $k-1,$ the height of the analysis
  lattice of $L''$ is at most $k-1$. Since the analysis lattice of $L$
  is obtained by grafting the analysis lattices of the vertex groups
  of the freely indecomposable free factors of $L$ (the groups $L''$)
  to the leaves at level $0.1$ in the analysis lattice of $L,$ its
  height is at most $k$. Since $k\leq 3n$ the theorem is proven.
\end{proof}

\section{Remarks}

\par There is an easier proof of Theorem~\ref{strictlengthbound} which
does not use modular group and generalizes to limit groups over linear
groups. One uses equational Noetherianness, the fact that the
representation variety has only finitely many components, and a
diagonal argument to achieve the same end. 

The preprint~\cite{houcine-2008} of Ould-Houcine is written in more
model theoretic language than this paper, and uses an argument like
the above. Theorem~\ref{strictlengthbound} follows from his bound on
the Cantor-Bendixon rank of the closure of the space of marked free
groups plus finite presentability. To clarify, if a sequence of marked
limit groups $(G,S_i)$ converges to a marked limit group $(H,S),$ then
$H$ has a strict homomorphism onto $G$. Conversely, let $S$ be a
generating set for $H,$ and fix a strict homomorphism $\pi\colon
H\onto G$. Suppose $f_i\colon G\to\free$ converges to $G$. If
$f_i\circ\pi\circ\phi_i,$ $\phi_i\in\Mod(H),$ converges to $H,$ then
the sequence of marked limit groups $(G,\pi\circ\phi_i(S))$ converges
to the marked limit group $(H,S)$.

\bibliographystyle{amsalpha} 
\bibliography{krull}

\end{document}